\documentclass[11pt]{article}
\usepackage{}
\usepackage{mathrsfs}
\usepackage{amssymb}
\usepackage{amsfonts}
\usepackage{color,xcolor}
\usepackage{graphicx}
\usepackage{manfnt}
\usepackage{graphicx}
\newcommand{\lbl}[1]{\label{#1}}
\newcommand{\be}{\begin{equation}}
\newcommand{\ee}{\end{equation}}
\newcommand\bes{\begin{eqnarray}}
\newcommand\ees{\end{eqnarray}}
\newcommand{\bess}{\begin{eqnarray*}}
\newcommand{\eess}{\end{eqnarray*}}
\newtheorem{theo}{\textbf{\ \ \quad Theorem}}[section]

\newtheorem{lem}{\textbf{\ \ \quad Lemma}}[section]

\newtheorem{prop}{\textbf{\ \ \quad Proposition}}[section]
\newtheorem{defi}{\textbf{\ \ \quad Definition}}[section]

\newtheorem{remark}{Remark}[section]

\newcommand{\nm}{\nonumber}

\def\[{{\Big[}}\def\]{{\Big]}}\def\<{{\langle}}\def\>{{\rangle}}\def\({{\Big(}}
\def\){{\Big)}}

\def\sgn{\mbox{\rm sgn}}

\def\={&\!\!=\!\!&}

\def\mR{{\mathbb R}}

\def\geq{\geqslant}\def\leq{\leqslant}
\def\div{\mathord{{\rm div}}}

\begin{document}
\title{\bf Kinetic solutions for \\ nonlocal stochastic conservation laws}
\author{Guangying Lv$^{a,b}$, Hongjun Gao$^a$ and Jinlong Wei$^c$\thanks{Corresponding author Email: weijinlong.hust@gmail.com} }
\date{$^a$ Institute of Mathematics, School of Mathematical Science\\
 Nanjing Normal University, Nanjing 210023, China\\
 $^b$ School of Mathematics and Statistics, Henan University\\
  Kaifeng, Henan 475001, China\\
$^c$ School of Statistics and Mathematics, Zhongnan University of
\\ Economics and Law, Wuhan, Hubei 430073, China }
 \maketitle
\noindent{\hrulefill}
\vskip1mm\noindent
{\bf Abstract} This work is devoted to examine the uniqueness and existence of kinetic solutions for
 a class of scalar conservation laws involving a nonlocal super-critical diffusion operator
 and a multiplicative noise. Our proof for uniqueness is based upon the analysis on
 double variables method and the existence is enabled by a parabolic approximation.


 \vskip0mm\noindent{\hrulefill}

\section{Introduction}\label{sec1}\setcounter{equation}{0}

The present paper is concerned with the  anomalous  diffusion related to
the L\'{e}vy flights \cite{Str,SZF, Duan2015}.
 At the macroscopic modeling level, this means the   Laplacian for normal diffusion is replaced by a fractional power of the (negative) Laplacian.  We consider the following  partial differential equation, coupling   a conservation law with an anomalous  diffusion:
   \bes
\partial_tu(x,t)+\nu(-\Delta_x)^{\frac{\alpha}{2}} u+\div(A(u))=\Phi(u)\partial_tW(t), \ \
x\in\mathbb{R}^d,\ t\in (0,T),
    \lbl{1.1}\ees
with initial data:
   \bes
u(0,x)=u_0(x), \ \ \ \ x\in\mathbb{R}^d,
   \lbl{1.2}\ees
where $\nu$ is a nonnegative parameter, $\alpha\in(0,1)$,
and $A=(A_1,\cdots,A_d)$, a vector field (\emph{the flux}), is supposed to be
of class $C^2$ and its derivatives have at most polynomial growth. Following \cite{DV},
we assume that
$W$ is a cylindrical Wiener process: $W=\sum_{k\geq1}\beta_ke_k$, where $\beta_k$
are independent Brownian process and $\{e_k\}_{k\geq1}$ is a complete orthonormal
basis in a Hilbert space $H$. For each $u\in\mathbb{R}$, $\Phi(u):H\to L^2(\mathbb{R}^d)$
is defined by $\Phi(u)e_k=g_k(\cdot,u)$, where $g_k(\cdot,u)$ is a regular function on $\mathbb{R}^d$.
More precisely, we assume $g_k\in C(\mathbb{R}^{d+1})$ with the bounds
    \bes
&&G^2(x,u)=\sum_{k\geq1}|g_k(x,u)|^2\leq D_0(\hat g(x)+|u|^2),\lbl{1.3}\\
&&\sum_{k\geq1}|g_k(x,u)-g_k(y,v)|^2\leq D_1(|x-y|^2+|u-v|h(|u-v|)),
 \lbl{1.4}\ees
where $0\leq \hat g(x)\in L^1(\mathbb{R}^d),\, x,y\in\mathbb{R}^d,\,u,v\in\mathbb{R}$,
and $h$ is a continuous non-decreasing  function
on $\mathbb{R}_+$ with $h(0)=0$.

We briefly mention some    recent works on well-posedness of (\ref{1.1})-(\ref{1.2}),
which are relevant for the present paper.
To add a stochastic forcing $\Phi(u)dw(t)$ is natural for applications, which appears
in wide variety of field as physics, engineering, biology and so on.
We first recall some results on the stochastic scalar conservation law without diffusion ($\nu=0$):
\begin{eqnarray}\label{1.5}
\partial_t u(x,t)+\mbox{div}_xA(u)=\Phi(u)\partial_tW(t), \ \ x\in \mR^d, \ t\in (0,T).
\end{eqnarray}
 The Cauchy problem
of equation (\ref{1.5}) with additive noise has been studied in \cite{K}, where J. U. Kim
proposed a method of compensated compactness to prove, via vanishing viscosity approximation,
the existence of a stochastic weak entropy solution. A Kruzhkov-type method was used to prove the
uniqueness. Vallet-Wittbold \cite{VW} extended the results of Kim to the multi-dimensional Dirichlet
problem with additive noise. By using vanishing viscosity method, Young measure techniques and
Kruzhkov doubling variables technique, they proved the existence and uniqueness of the stochastic
entropy solution.

Concerning multiplicative noise, for Cauchy problem, Feng-Nualart \cite{FN} introduced a notion of
strong entropy solution in order to prove the uniqueness for the entropy solution of (\ref{1.5}). Using the vanishing
viscosity and compensated compactness arguments, they established the existence of stochastic
strong entropy solution only in $1D$ case. Chen et al. \cite{CDK} proved that the multi-dimensional
stochastic problem is well-posedness by using a uniform spatial BV-bound. Following the idea of
\cite{FN,CDK}, Lv et al. \cite{LDG} considered the Cauchy problem (\ref{1.1}).
Bauzet et al.\cite{BVW} proved a result of existence and uniqueness of the weak measure-valued
entropy solution to the multi-dimensional Cauchy problem (\ref{1.5}).

Using a kinetic formulation, Debussche-Vovelle \cite{DV} obtained a result of existence and
uniqueness of the entropy solution to the problem posed in a d-dimensional torus.
About the Cauchy-Dirchlet problem (\ref{1.5}), see \cite{BVW2014}.

When $\nu>0$, the problem (\ref{1.1}) with (\ref{1.2}) has been studied in
\cite{LDG,LDGW}, where the Kruzhkov's semi-entropy formulations was used.
It is remarked that except for \cite{DV,FN}, the previous results only considered
the Brownian motion perturbation. That is, the noise does not depend on the
spatial variable.

Inspired by \cite{DV}, in this paper, we reconsider the problem (\ref{1.1}) with
(\ref{1.2}) and obtain the well-posedness by using the kinetic formulation.
The advantage of kinetic formulation method is that we can deal with the
cylindrical Wiener process in any dimension.
We remark that it is not trivial to generalize the results of \cite{DV} to the problem (\ref{1.1})
with (\ref{1.2}). Because of the nonlocal term $(-\Delta_x)^{\frac{\alpha}{2}} u$, the proof of existence for kinetic solutions will become more complicated and the assumptions
on the initial data will become stronger. Moreover, compared with \cite{DV}, we have to introduce another non-negative measure to overcome the difficulty.
The proof of uniqueness of solution to (\ref{1.1}) will be
different from that in (\ref{1.5}). In this paper, we mainly focus on
how to deal with the nonlocal term.

This paper is organized as follows. In Section
2, we introduce some notions on solutions for (\ref{1.1})-(\ref{1.2}),
and then prove the uniqueness and existence of  kinetic solutions in Section 3.
We further discuss
the regularity properties and continuous dependence (on nonlinearities and L\'{e}vy measures)
for kinetic solutions in Section 4.

\section{Entropy solutions and kinetic solutions}\label{sec2}
\setcounter{equation}{0}
In this section, we first give the definitions of stochastic entropy solutions and
stochastic kinetic solutions, then prove that they are equivalent and last state out
our main results.

To present our formulation for (\ref{1.1}), we recall the following results
on the operator $(-\Delta)^{\frac{\alpha}{2}}$.

\begin{lem}\lbl{l2.1} {\rm(\cite{DI})} There exists a constant $C_d(\alpha)>0$
that only depends on $d$ and $\alpha$, and such that for all $\phi\in\mathcal {S}(\mathbb{R}^d)$,
all $r>0$ and all $x\in\mathbb{R}^d$
   \bess
(-\Delta)^{\frac{\alpha}{2}}\phi(x)&=&-C_d(\alpha)\int_{|z|\geq r}\frac{\phi(x+z)-\phi(x)}{|z|^{d+\alpha}}dz\nm\\
&&
-C_d(\alpha)\int_{|z|\leq r}\frac{\phi(x+z)-\phi(x)-\nabla\phi(x)\cdot z}{|z|^{d+\alpha}}dz.
    \eess
Moreover, when $\alpha\in(0,1)$, one can take $r=0$.
\end{lem}

We take $\nu=1$ in Sections 2 and 3.
Here and in the followings, we use $(\cdot,\cdot)$ to denote the
inner product of $L^2$-valued functions. Following \cite{FN}, we have the definition.

\begin{defi}\lbl{d2.1} {\rm(Stochastic Nonlocal Entropy solution)} An $L^2(\mathbb{R}^d)$-valued $\{\mathcal {F}_t:0\leq t\leq T\}$-predictable stochastic process $\{u(t)=u(x,t)\}$ is called a
stochastic entropy solution of {\rm(\ref{1.1})} provided

{\rm(1)} For each $p\geq 1$
   \bess
\mathbb{E}[ess\sup_{0\leq t\leq T}\|u(t)\|_p^p]<\infty;
  \eess

{\rm(2)} For $0\leq\psi\in C_c^{1,2}([0,T)\times\mathbb{R}^d)$ and all
convex $\eta\in C^2(\mathbb{R})$, the following inequality holds
\bes
&&\int^T_0(\eta(u(r)),\partial_t\psi(r,\cdot)) dr+(\eta(u(0)),\psi(0,\cdot))+\int^T_0(\Psi(u(r)),\nabla_x\psi(r,\cdot)) dr\nm\\
&&
-\int^T_0(\eta(u(r)),(-\Delta)^{\frac{\alpha}{2}}\psi(r)) dr
+\sum_{k\geq1}\int^T_0(g_k(\cdot,u(r))\eta'(u(r)),\psi(r,\cdot)) dW(r)\nm\\
&&+\frac{1}{2}\int^T_0( G^2(\cdot,u(r))\eta''(u(r,\cdot)),\psi(r,\cdot)) dr\geq 0,
   \lbl{2.1}\ees
a.s., where $\Psi(u(r))=\int_0^ua(\xi)\eta'(\xi)d\xi$ and $a(\xi)=A'(\xi)$.
\end{defi}

\begin{remark}\lbl{r2.1} Comparing the above Definition \ref{d2.1} with
Definition 2.1 in \cite{A}, we get that the above Definition \ref{d2.1} is the case
"intermediate".
It is well-known that "classical $\Rightarrow$ entropy $\Rightarrow$ intermediate $\Rightarrow$ weak",
see Remark 4.2 in \cite{A}.

We cannot define the solution as the Definition 2.1 in \cite{A},
because we cannot prove the solutions of (\ref{1.1}) belong to $BV(\mathbb{R}^d)$
even if the initial data belong to $BV(\mathbb{R}^d)$. When $\nu=0$ in (\ref{1.1}),
Chen et al. \cite{CDK} obtain the Fractional BV-estimate, see \cite[Theorem 7]{CDK}.
Thus it is impossible to get the BV-estimate of solution to (\ref{1.1}) under the assumption
that the noise term depends on the spatial variable $x$.
Besides, the BV-estimate of solution to (\ref{1.1}) with $\Phi\equiv0$ was obtained in
\cite{WDL}. That is, the deterministic nonlocal conservation law keeps the Bounded Variation property.
\end{remark}

\begin{remark}\lbl{r2.2} The solution defined in Definition \ref{d2.1} satisfies
the initial condition in the following sense: for any compact set $K\subset\mathbb{R}^d$,
   \bess
ess\lim\limits_{t\to0+}\mathbb{E}\int_K|u-u_0|dx=0.
  \eess
The proof is exactly as that of Remark 2.7 in \cite{BVW}.
 \end{remark}

Inspired by \cite{DV}, we give the following definitions.
  \begin{defi}\lbl{d2.2} {\rm(Kinetic \,measure)} We say
that a map $m$ from $\Omega$ to the set of non-negative finite measure
over $\mathbb{R}^d\times[0,T]\times\mathbb{R}$ is a kinetic measure if

1. $m$ is measurable, in the sense that for each
$\phi\in C_b(\mathbb{R}^d\times[0,T]\times\mathbb{R})$,
$\langle m, \phi\rangle:\Omega\to\mathbb{R}$ is,

2. $m$ vanishes for large $\xi$: if $B_R^c=\{\xi\in\mathbb{R},|\xi|\geq R\}$, then
  \bess
\lim\limits_{R\to\infty}\mathbb{E}m(\mathbb{R}^d\times[0,T]\times B_R^c)=0,
   \eess

3. for all $\phi\in C_b(\mathbb{R}^{d+1})$, the process
   \bess
t\mapsto\int_{\mathbb{R}^d\times[0,t]\times\mathbb{R}}\phi(x,\xi)dm(x,s,\xi)
   \eess
is predictable.
 \end{defi}

\begin{defi}\lbl{d2.3}{\rm(Solution)} Let
$u_0\in L^\infty\cap  L^1(\mathbb{R}^d)$. A measurable function
$u:\mathbb{R}^d\times[0,T]\times\Omega\to\mathbb{R}$ is said to be a solution to
(\ref{1.1}) with initial datum $u_0$ if $\{u(t)\}$ is predictable, for all $p\geq1$,
there exists $C_p\geq0$ such that
   \bess
\mathbb{E}\left(ess \sup_{t\in[0,T]}\|u(t)\|_{L^p(\mathbb{R}^d)}^p\right)\leq C_p,
  \eess
and if there exists a kinetic measure $m$ such that $f:=\textbf{1}_{u>\xi}$ satisfies:
for all $\varphi\in C_c^2(\mathbb{R}^d\times[0,T)\times\mathbb{R})$,
    \bes
&&\int_0^T\langle f(t),\partial_t\varphi(t)\rangle dt+\langle f_0,\varphi(0)\rangle
+\int_0^T\langle f(t),a(\xi)\cdot\nabla\varphi(t)\rangle dt\nm\\
&=&\int_0^T\langle f(t),(-\Delta_x)^{\frac{\alpha}{2}}\varphi(t)\rangle dt
-\sum_{k\geq1}\int_0^T\int_{\mathbb{R}^d}g_k(x,u(x,t))\varphi(x,t,u(x,t))dxd\beta_k(t)\nm\\
&&-\frac{1}{2}\int_0^T\int_{\mathbb{R}^d}\partial_\xi\varphi(x,t,u(x,t))G^2(x,u(x,t))dxdt+
m(\partial_\xi\varphi),
   \lbl{2.2}\ees
a.s., where $f_0(x,\xi)=\textbf{1}_{u_0(x)>\xi}$, $G^2:=\sum_{k=1}^\infty|g_k|^2$ and $a(\xi)=A'(\xi)$.
\end{defi}

In (\ref{2.2}), we have used the brackets $\langle\cdot,\cdot\rangle$ to denote
the duality between $C_c^\infty(\mathbb{R}^{d+1})$ and the space of distributions
over $\mathbb{R}^{d+1}$. In what follows, we will denote similarly the integral
   \bess
\langle F,G\rangle=\int_{\mathbb{R}^{d+1}}F(x,\xi)G(x,\xi)dxd\xi,\ \ F\in L^p(\mathbb{R}^{d+1}),\,
G\in L^q(\mathbb{R}^{d+1}),
   \eess
where $1\leq p<\infty$ and $q$ is the conjugate exponent of $p$. In (\ref{2.2}),
we also have indicated the dependence of $g_k$ and $G^2$ on $u$, which is actually absent
in the additive case and we have used the shorthand $m(\psi)$ for
   \bess
m(\psi)=\int_{\mathbb{R}^d\times[0,T]\times\mR}\psi(x,t,\xi)dm(x,t,\xi),\ \ \ \psi\in C_c(\mathbb{R}^d\times[0,T]\times\mR).
   \eess
Equation (\ref{2.2}) is the weak form of the equation
    \bes
\left(\partial_t+a(\xi)\cdot\nabla_x+(-\Delta_x)^{\frac{\alpha}{2}}\right)\textbf{1}_{u(x,t)>\xi}
=\delta_{u=\xi}\Phi\dot{W}+\partial_\xi(m-\frac{1}{2}G^2\delta_{u=\xi}).
   \lbl{2.3}\ees
Now, we present a formal derivation of equation (\ref{2.3}) from (\ref{1.1}) for regular solution which is similar to \cite{DV}. It is essentially a consequence of It\^{o} formula.
Indeed, by the identity
$(\textbf{1}_{u>\xi},\theta'):=\int_\mathbb{R}\textbf{1}_{u>\xi}\theta'd\xi=\theta(u)-\theta(-\infty)$
for $\theta\in C^\infty(\mathbb{R})$, it yields that
   \bess
d(\textbf{1}_{u>\xi},\theta')&=&\theta'(u)\left(-a(u)\cdot\nabla u-(-\Delta_x)^{\frac{\alpha}{2}}u
+\Phi(u)dW\right)+\frac{1}{2}\theta''(u)G^2dt\\
&=&-\div\left(\int^ua(\xi)\theta'(\xi)d\xi\right)dt+\frac{1}{2}\theta''(u)G^2dt+\theta'(u)\Phi(u)dW\\
&&-
(-\Delta_x)^{\frac{\alpha}{2}}\theta(u)
-\theta'(u)(-\Delta_x)^{\frac{\alpha}{2}}u+(-\Delta_x)^{\frac{\alpha}{2}}\theta(u)\\
&=&-\div(a\textbf{1}_{u>\xi},\theta')dt-\frac{1}{2}(\partial_\xi(G^2\delta_{u=\xi}),\theta')
+(\delta_{u=\xi}\Phi dW,\theta')\\
&&-((-\Delta_x)^{\frac{\alpha}{2}}\textbf{1}_{u>\xi},\theta')+(\partial_\xi m,\theta'),
   \eess
where $(\partial_\xi m,\theta')=(-\Delta_x)^{\frac{\alpha}{2}}\theta(u)-\theta'(u)(-\Delta_x)^{\frac{\alpha}{2}}u$,
 we have used the following fact
    \bess
(-\Delta_x)^{\frac{\alpha}{2}}\theta(u)&=&
C_d(\alpha)\int_{\mathbb{R}^d\setminus\{0\}}\frac{\theta(u)-\theta(u(x+z,t))}{|z|^{d+\alpha}}dz\\
&=&C_d(\alpha)\int_{\mathbb{R}^d\setminus\{0\}}
\frac{1}{|z|^{d+\alpha}}\left(\int_{u(x+z,t)}^{u}\theta'(\xi)d\xi\right)dz\\
&=&C_d(\alpha)\int_{\mathbb{R}^d\setminus\{0\}}
\frac{(\textbf{1}_{u>\xi}-\textbf{1}_{u(x+z,t)>\xi},\theta'(\xi))}{|z|^{d+\alpha}}dz\\
&=&((-\Delta_x)^{\frac{\alpha}{2}}\textbf{1}_{u>\xi},\theta').
   \eess
Next, we calculate $m$. Firstly, let $\theta\in C_c^2(\mathbb{R})$ be a convex function and we have
  \bess
\langle\partial_\xi m,\theta'\rangle=-\langle m,\theta''\rangle=-\int_\mathbb{R}m(x,t,\xi)\theta''(\xi)d\xi.
   \eess
Assume that $\theta_\epsilon\in C_c^2(\mathbb{R})$  is a convex function satisfying
 $\lim\limits_{\epsilon\to0}\theta_\epsilon(\xi)=|v-\xi|$,
 $\lim\limits_{\epsilon\to0}\theta_\epsilon'(\xi)=\sgn(\xi-v)$
 and $\lim\limits_{\epsilon\to0}\theta_\epsilon''(\xi)=\delta_{v=\xi}$,
 where $v\in\mathbb{R}$, then we get
   \bess
m(x,t,v)&=&-\lim\limits_{\epsilon\to0}\langle\partial_\xi m,\theta_\epsilon'\rangle\\
&=&-\lim\limits_{\epsilon\to0}
(-\Delta_x)^{\frac{\alpha}{2}}\theta_\epsilon(u)-\theta_\epsilon'(u)(-\Delta_x)^{\frac{\alpha}{2}}u\\
&=&\sgn(u-v)(-\Delta_x)^{\frac{\alpha}{2}}u-(-\Delta_x)^{\frac{\alpha}{2}}|u-v|.
  \eess
From the definition of $(-\Delta)^{\frac{\alpha}{2}}$, we have another representation for the kinetic measure $m$. It follows from Lemma \ref{l2.1} that
   \bess
(\partial_\xi m,\theta')&=&(-\Delta_x)^{\frac{\alpha}{2}}\theta(u)-\theta'(u)(-\Delta_x)^{\frac{\alpha}{2}}u\\
&=&C_d(\alpha)\int_{\mathbb{R}^d\setminus\{0\}}
\frac{1}{|z|^{d+\alpha}}[\theta(u)-\theta(u(x+z,t))-\theta'(u)(u-u(x+z,t))]dz\\
&=&\frac{C_d(\alpha)}{2}\int_{\mathbb{R}^d\setminus\{0\}}
\frac{1}{|z|^{d+\alpha}}\theta''((1-\tau)u+\tau u(x+z,t))(u-u(x+z,t))^2dz\\
&=&-\frac{C_d(\alpha)}{2}\int_{\mathbb{R}^d\setminus\{0\}}
\frac{1}{|z|^{d+\alpha}}\left(\partial_\xi((u-u(x+z,t))^2
\delta_{\xi=(1-\tau)u+\tau u(x+z,t)},\theta'(\xi)\right)dz,
  \eess
which implies that
    \bes
m(x,t,\xi)=\frac{C_d(\alpha)}{2}\int_{\mathbb{R}^d\setminus\{0\}}
\frac{(u-u(x,t+z))^2}{|z|^{d+\alpha}}\delta_{\xi=(1-\tau)u+\tau u(x+z,t)}dz,\ \ \ \tau\in(0,1)
   \lbl{2.4}\ees
The above representation shows that $m_1$ is a non-negative measure.

Taking $\theta(\xi)=\int_{-\infty}^\xi\varphi$,
we then obtain the formulation. Besides the nonnegative measure given by (\ref{2.4}), the kinetic measure $m$ described in (\ref{2.2}) contains another non-negative measure, which is sometimes interpreted as a Lagrange
multiplier for the evolution of $f$ by $\partial_t+a\cdot\nabla$ under the constraint
$f={\rm graph}=\textbf{1}_{u>\xi}$. It will be arose when $u$ becomes discontinuous. Indeed, if one we add the viscosity term $\varepsilon\Delta u$ in equation (\ref{1.1}), then the measure $m$ can be written as
    \bess
m(\phi)&=&\varepsilon\int_{\mathbb{R}^{d}\times[0,T]}\phi(x,t,u(x,t))|\nabla u|^2dxdt\cr\cr&&+\int_{\mathbb{R}^{d}\times[0,T]\times\mR}\phi(x,t,\xi)\frac{C_d(\alpha)}{2}
\int_{\mathbb{R}^d\setminus\{0\}}
\frac{(u-u(x,t+z))^2}{|z|^{d+\alpha}}\delta_{\xi=(1-\tau)u+\tau u(x+z,t)}dzdxdtd\xi.
   \eess

Now, we are in a position to show the relationship between entropy
solutions and kinetic solutions for (\ref{1.1})- (\ref{1.2}).

\begin{theo} \lbl{t2.1}{\rm(Kinetic formulation)} Let $u_0\in L^1\cap L^\infty(\mR^d)$.
For a measurable function $u:\mathbb{R}^d\times[0,T]\times\Omega\to\mathbb{R}$, it
is equivalent to be a kinetic solution to (\ref{1.1}), i.e. both the solutions in sense of
Definitions \ref{d2.1} and \ref{d2.3} are equivalent.
\end{theo}

\textbf{Proof.} Choosing test function $\varphi(x,t,\xi)=\psi(x,t)\eta'(\xi)$ in
(\ref{2.2}) and noting that $\eta$ is a convex function, we have
    \bess
\langle f(t),(-\Delta_x)^{\frac{\alpha}{2}}\varphi(t)\rangle
&=&\langle(-\Delta_x)^{\frac{\alpha}{2}} f(t),\psi(x,t)\eta'(\xi)\rangle\\
&=&((-\Delta_x)^{\frac{\alpha}{2}} \eta(u(t)),\psi(x,t))\\
&=&( \eta(u(t)),(-\Delta_x)^{\frac{\alpha}{2}}\psi(x,t)).
  \eess
Using the above inequality and the facts $m(\eta'')\geq0$ and $n(\eta'')\geq0$,
(\ref{2.2}) implies the inequality in Definition \ref{d2.1}. That is,
a kinetic solution will be a entropy solution.

Conversely, similar to \cite{DV}, one defines the measure $m$ by
    \bess
m(\eta''\psi)&=&\int_0^T(\eta(u),\partial_t\psi)dr+(\eta(u_0),\psi(0))
+\int_0^T(\Psi(u),\nabla\psi)dr\\
&&
+\sum_{k\geq1}\int_0^T(g_k(\cdot,u(r))\eta'(u(r)),\psi)d\beta_k(r)
+\frac{1}{2}\int_0^T(G^2(\cdot,u(r))\eta''(u(r)),\psi)dr\\
&&-\int_0^T\int_{\mR^d}\eta(u(r))(-\Delta)^\frac\alpha2\psi dxdr,
   \eess
then one derives (\ref{2.1}). Moreover, by virtue of above representation of $m$, we prove that $m$ is a kinetic measure. $\Box$

In order to prove the existence of solution, we introduce the following definitions, see \cite{DV}.

\begin{defi} \lbl{d2.4}{\rm(Young measure)} Let $(X,\lambda)$ be a finite measure space.
Let $\mathcal {P}_1(\mathbb{R})$ denote the set of probability measures on $\mathbb{R}$. We
say that a map $\nu:X\to\mathcal {P}_1(\mathbb{R})$ is a Young measure on $X$ if, for all
$\phi\in C_b(\mathbb{R})$, the map $z\mapsto\nu_z(\phi)$ from $X$ to $\mathbb{R}$ is measurable.
We say that a Young measure $\nu$ vanishes at infinity if, for every $p\geq1$,
   \bes
\int_X\int_\mathbb{R}|\xi|^pd\nu_z(\xi)d\lambda(z)<+\infty.
   \lbl{2.5}\ees
\end{defi}

\begin{defi}\lbl{d2.5} {\rm(Kinetic function)}Let $(X,\lambda)$ be a finite measure space.
A measurable function $f:X\times\mathbb{R}\to[0,1]$ is said to be a kinetic function if there
exists a Young measure $\nu$ on $X$ that vanishes at infinity such that, for
$\lambda-a.e.\,z\in X$, for all $\xi\in\mathbb{R}$,
  \bess
f(z,\xi)=\nu_z(\xi,+\infty).
  \eess
We say that $f$ is an equilibrium if there exists a measurable function $u:\,X\to\mathbb{R}$
such that $f(z,\xi):=\textbf{1}_{u(z)>\xi}$ a.e., or, equivalently, $\nu_z=\delta_{u(z)}$ for
a.e. $z\in X$.
  \end{defi}

If $f:X\times\mathbb{R}\to[0,1]$ is a kinetic function, we denote by $\bar f$ the
conjugate function $\bar f=1-f$. We can define the kinetic function in another way (see \cite{LPT})
    \bess
\chi_u(\xi)&=&\textbf{1}_{(0,u(x,t))}(\xi)-\textbf{1}_{(u(x,t),0)}(\xi)=\textbf{1}_{u>\xi}-\textbf{1}_{0>\xi},
   \eess
which is decreasing faster than any power of $\xi$ at infinity. Contrary to $f$,
$\chi_u(\xi)$ is integrable. Now, we recall the compactness of Young measures,
see \cite{DV} for the proof.
   \begin{prop}\lbl{p2.1} {\rm\cite[Theorem 5]{DV}} Let $(X,\lambda)$ be a finite measure space
such that $L^1(X)$ is separable. Let $(\nu^n)$ be a sequence of Young measures on
$X$ satisfying (\ref{2.5}) uniformly for some $p\geq1$:
    \bes
\sup_n\int_X\int_\mathbb{R}|\xi|^pd\nu^n_z(\xi)d\lambda(z)<+\infty.
   \lbl{2.6}\ees
Then there exists a Young measure $\nu$ on $X$ and a subsequence still denoted
$(\nu^n)$ such that, for all $h\in L^1(X)$, for all $\phi\in C_b(\mathbb{R})$,
   \bess
\lim\limits_{n\to\infty}\int_Xh(z)\int_{\mathbb{R}}\phi(\xi)d\nu_z^nd\lambda(z)
=\int_Xh(z)\int_{\mathbb{R}}\phi(\xi)d\nu_zd\lambda(z).
   \eess
\end{prop}

By Proposition \ref{p2.1}, we have the following result. Let $(f_n)$ be a
sequence of kinetic functions on $X\times\mathbb{R}$: $f_n(z,\xi)=\nu_z^n(\xi,+\infty)$,
where $\nu^n$ are Young measures on $X$ satisfying (\ref{2.6}). Let $f$ be a kinetic function
on $X\times\mathbb{R}$ such that $f_n\rightharpoonup f$ in $L^\infty(X\times\mathbb{R})$
weak$-^*$. Assume that $f_n$ and $f$ are equilibria:
   \bess
f_n(z,\xi)=\textbf{1}_{u_n(z)>\xi}, \ \ \ \ f(z,\xi)=\textbf{1}_{u(z)>\xi}.
  \eess
Then, for all $1\leq q<p$, $u_n\to u$ in $L^q(X)$ strong.

\begin{defi}\lbl{d2.6} {\rm(Generalized solution)} Let $f_0:\Omega\times\mathbb{R}^{d+1}\mapsto[0,1]$
be a kinetic function. A measurable function $f:\Omega\times\mR^d\times[0,T]\times\mR\mapsto[0,1]$
is said to be a generalized solution to (\ref{1.1}) with initial datum $f_0$ if $\{f(t)\}$
is predictable and is a kinetic function such that: for all $p\geq1$, $\nu:=-\partial_\xi f$
satisfies
   \bes
\mathbb{E}\left(ess\sup_{t\in[0,T]}\int_{\mathbb{R}^{d+1}}|\xi|^pd\nu_{x,t}dx\right)\leq C_p,
   \lbl{2.7}\ees
where $C_p$ is a positive constant and: there exists a kinetic measure $m$ such that
for all $\varphi\in C_c^2(\mR^d\times[0,T)\times\mR)$,
   \bes
&&\int_0^T\langle f(t),\partial_t\varphi(t)\rangle dt+\langle f_0,\varphi(0)\rangle
+\int_0^T\langle f(t),a(\xi)\cdot\nabla\varphi(t)\rangle dt\nm\\
&=&\int_0^T\langle f(t),(-\Delta_x)^{\frac{\alpha}{2}}\varphi(t)\rangle dt
-\sum_{k\geq1}\int_0^T\int_{\mathbb{R}^{d+1}}g_k(x,\xi)\varphi(x,t,\xi)d\nu_{x,t}(\xi)dxd\beta_k(t)\nm\\
&&-\frac{1}{2}\int_0^T\int_{\mathbb{R}^{d+1}}\partial_\xi\varphi(x,t,\xi)G^2(x,\xi)d\nu_{x,t}(\xi)dxdt
+m(\partial_\xi\varphi), \ \ a.s..
   \lbl{2.8}\ees
\end{defi}

Note that the generalized solution (Definition \ref{d2.6}) implies the solution
(Definition \ref{d2.3}). Indeed, if $f$ is a generalized solution such that
$f=\textbf{1}_{u>\xi}$, then $u(x,t)=\int_\mathbb{R}(f-\textbf{1}_{0>\xi})d\xi$,
hence $u$ is predictable. Moreover, if $\nu_{x,t}(\xi)=\delta_{u=\xi}$, then equality
(\ref{2.8}) implies (\ref{2.2}).

Following \cite{DV}, we shall show that any generalized solution admits
possibly different left and right weak limits at any point $t\in[0,T]$ almost
surely. This property is important to prove a comparison principle which
allows to prove uniqueness. Meanwhile, it allows us to rewrite (\ref{2.8})
in some stronger sense.
   \begin{prop}\lbl{p2.2}{\rm(left and right weak limits)} Let $f_0$ be a
kinetic initial datum. Let $f$ be a generalized solution to (\ref{1.1}) with
initial datum $f_0$. Then $f$ admits almost surely left and right limits at
all points $t_*\in[0,T]$. More precisely, for all $t_*\in[0,T]$ there exists some
kinetic functions $f^{*,\pm}$ on $\Omega\times\mathbb{R}^{d+1}$ such that $\mathbb{P}$-a.s.
    \bess
\langle f(t_*-\varepsilon),\varphi\rangle \to\langle f^{*,-},\varphi\rangle
    \eess
and
    \bess
\langle f(t_*+\varepsilon),\varphi\rangle \to\langle f^{*,+},\varphi\rangle
    \eess
as $\varepsilon\to0$ for all $\varphi\in C_c^2(\mathbb{R}^{d+1})$. Moreover, almost
surely,

     \bes
\langle f^{*,+}-f^{*,-},\varphi\rangle=-\int_{\mR^d\times[0,T]\times\mR}\partial_\xi\varphi(x,\xi)
\textbf{1}_{\{t_*\}}(t)dm(x,t,\xi).
     \lbl{2.9} \ees
In particular, almost surely, the set of $t_*\in[0,T]$ such that $f^{*,+}\neq f^{*,-}$ is
countable.
   \end{prop}

{\bf Proof.} Following \cite{DV}, for all $\varphi\in C_c^2(\mathbb{R}^{d+1})$, a.s.,
the map
    \bess
J_\varphi:t&\to&\int_0^t\langle f(s),a(\xi)\cdot\nabla\varphi \rangle ds
-\int_0^t\langle f(s),(-\Delta_x)^{\frac{\alpha}{2}}\varphi \rangle ds\nm\\
&&\qquad
+\sum_{k\geq1}\int_0^t\int_{\mathbb{R}^{d+1}}g_k(x,\xi)\varphi(x,\xi)d\nu_{x,s}(\xi)dxd\beta_k(s)\nm\\
&&\qquad\qquad+\frac{1}{2}\int_0^t\int_{\mathbb{R}^{d+1}}\partial_\xi\varphi(x,\xi)G^2(x,\xi)d\nu_{x,s}(\xi)dxds
   \eess
is continuous on $[0,T]$. Taking the test function of the form $(x,t,\xi)\mapsto\varphi(x,\xi)\gamma(t)$,
$\gamma\in C_c^1([0,T]),\,\varphi\in C_c^2(\mathbb{R}^{d+1})$, we get by using Fubini Theorem
and the weak formulation (\ref{2.8})
   \bess
\int_0^T\mathcal {J}_\varphi(t)\gamma'(t)dt+\langle f_0,\varphi\rangle\gamma(0)=\langle m,\partial_\xi\varphi\rangle(\gamma),
   \eess
where $\mathcal {J}_\varphi(t):=\langle f(t),\varphi\rangle-J_\varphi(t)$. This
shows that $\partial_t\mathcal {J}_\varphi$ is a measure on $(0,T)$, i.e.,
the function $\mathcal {J}_\varphi\in BV(0,T)$. Hence it admits left and
right limits at all points $t_*\in[0,T]$. Since $J_\varphi$ is continuous, this
also holds for $\langle f(t),\varphi\rangle$: for all $t_*\in[0,T]$, the limits
   \bess
\langle f,\varphi\rangle(t_*+):=\lim\limits_{t\downarrow t_*}\langle f,\varphi\rangle(t)\ \ {\rm and}\  \
\langle f,\varphi\rangle(t_*-):=\lim\limits_{t\uparrow t_*}\langle f,\varphi\rangle(t)
    \eess
exist. Then following the proof of Proposition 8 of \cite{DV}, it is easy to complete
the proof. $\Box$

Using the Proposition \ref{p2.2}, we can derive a kinetic formulation at given
$t$. Taking a test function of the form $(x,t,\xi)\mapsto [K(T-t,\cdot)\ast\varphi(\cdot,\xi)](x)\gamma(t)$
where the kernel function $K$ satisfies $K_t+(-\Delta)^{\frac{\alpha}{2}}K=0$, and $\gamma$ is the function
   \bess
\gamma(s)=\left\{\begin{array}{lllll}
1,\ \ \ & s\leq t,\\
1-\frac{s-t}{\varepsilon},\ \ \ &t\leq s\leq t+\varepsilon,\\
0, \ \ \ \ &t+\varepsilon\leq s,
\end{array}\right.\eess
we obtain at the limit $[\varepsilon\to0]:$ for all $t\in[0,T]$ and
$\varphi\in C_c^2(\mathbb{R}^{d+1})$,
      \bes
&&-\langle f^+(t),\tilde\varphi(t)\rangle +\langle f_0,\tilde\varphi(0)\rangle
+\int_0^t\langle f(s),a(\xi)\cdot\nabla\tilde\varphi \rangle ds\nm\\
&=&
-\sum_{k\geq1}\int_0^t\int_{\mathbb{R}^{d+1}}g_k(x,\xi)\tilde\varphi(x,\xi)d\nu_{x,s}(\xi)dxd\beta_k(s)\nm\\
&&-\frac{1}{2}\int_0^t\int_{\mathbb{R}^{d+1}}\partial_\xi\tilde\varphi(x,\xi)G^2(x,\xi)d\nu_{x,s}(\xi)dxds+
\langle m,\partial_\xi\tilde\varphi\rangle([0,t]), \ \ a.s.,
   \lbl{2.10}\ees
where
    \bess
\tilde\varphi(x,t)=[K(T-t,\cdot)\ast\varphi(\cdot,\xi)](x),\
\langle m,\partial_\xi\tilde\varphi\rangle([0,t])=\int_{\mathbb{R}^d\times[0,t]\times\mathbb{R}}
\partial_\xi\tilde\varphi(x,\xi)dm(x,s,\xi).
   \eess
We remark that if $\varphi\in C_c^2(\mathbb{R}^{d+1})$, then $(-\Delta)^{\frac{\alpha}{2}}\varphi$ makes
sense and also for $(-\Delta)^{\frac{\alpha}{2}}\tilde\varphi$. Using the following fact
     \bess
\langle f^+(t),\tilde\varphi(t)\rangle&=&\int_{\mathbb{R}^{d+1}}f^+(x,t,\xi)\tilde\varphi(x,t,\xi)dxd\xi\\
&=&\int_{\mathbb{R}^{d+1}}f^+(x,t,\xi)\int_{\mathbb{R}^d}K(T-t,y)\varphi(x-y,\xi)dydxd\xi\\
&=&\int_{\mathbb{R}^{d+1}}\int_{\mathbb{R}^d}K(T-t,y)f^+(x,t-y,\xi)dy\varphi(x,\xi)dxd\xi\\
&=&\langle\tilde f^+(t),\varphi\rangle,
    \eess
we can rewrite the equality (\ref{2.10}), that is, the convolution of $K$ and $\varphi$
can be changed into the convolution of $K$ and another function.

\begin{remark}\lbl{r2.3} (The case of equilibrium) Suppose that $f^{*,-}$ is at equilibrium in (\ref{2.9}): there is a random variable $u^*\in L^1(\mR^d)$ so that $f^{*,-}=\textbf{1}_{u*>\xi}$ a.s. Let $m^*$ denote the restriction of $m$ to $\mR^d\times\{t_*\}\times\mR$. We thus have
\bess
f^{*,+}-1_{u*>\xi}=\partial_\xi m^*.
\eess
By the condition 2 in Definition 2.2, one achieves that
\bess
\int_{\mR}(f^{*,+}(x,\xi)-\textbf{1}_{0>\xi})d\xi=
\int_{\mR}(\textbf{1}_{u*>\xi}-\textbf{1}_{0>\xi})d\xi=u^*.
\eess
Observing that
\bess
p^*: \xi\rightarrow \int_{-\infty}^{\xi}(\textbf{1}_{u*>\zeta}-f^{*,+}(\zeta))
\eess
is non-negative and $\partial_\xi(m^*+p^*)=0$, thus $m^*+p^*$ is constant and actually vanishes by the condition 2 in Definition 2.2 and the obvious fact that $p^*$ also vanishes when $\xi\rightarrow \infty$. Since $m^*\geq 0$, we conclude $m^*=0$, which suggests $f^{*,+}=f^{*,-}$.
 \end{remark}

\section{Uniqueness and existence of kinetic solutions}\label{sec3}\setcounter{equation}{0}

In this section, we are interested in the Cauchy problem (\ref{1.1})-(\ref{1.2}) and it is ready for us to state our main result.
\begin{theo}\lbl{t3.1}
Let (\ref{1.3}) and (\ref{1.4}) hold. Then there is a unique kinetic solution of the nonlocal Cauchy problem (\ref{1.1})-(\ref{1.2}).
\end{theo}

We will use the doubling variables method to prove the uniqueness. Let $f_i$ ($i=1,2$) be the generalized
solution of the equation
     \bes
\partial_tu_i(x,t)+(-\Delta_x)^{\frac{\alpha}{2}} u_i+\div(A(u_i))=\Phi(u_i)\partial_tW(t), \ \
x\in\mathbb{R}^d,\ t>0.
    \lbl{3.1}\ees
Set $\bar f=1-f$. In order to prove Theorem \ref{t3.1}, we need the following lemma.
   \begin{lem}\lbl{l3.1} Let $f_i,\,i=1,2$, be generalized solution
to (\ref{3.1}). Then, for $0\leq t\leq T$, and non-negative functions
$\rho\in C_c^\infty(\mathbb{R}^d)$ and $\psi\in C_c^\infty(\mathbb{R})$, we have
   \bes
&&\mathbb{E}\int_{\mathbb{R}^{2d}}\int_{\mathbb{R}^2}\breve\rho(x-y)\psi(\xi-\zeta)f_1^{\pm}(x,t,\xi)
\bar f_2^\pm(y,t,\zeta)d\xi d\zeta dxdy\nm\\
&\leq&\mathbb{E}\int_{\mathbb{R}^{2d}}\int_{\mathbb{R}^2}\breve\rho(x-y)\psi(\xi-\zeta)f_{1,0}^{\pm}(x,\xi)
\bar f_{2,0}^\pm(y,\zeta)d\xi d\zeta dxdy+I_{\breve\rho}+I_\psi,
   \lbl{3.2}\ees
where ,
    \bess
&&\quad\quad\breve\rho(x)=\int_{\mathbb{R}^d}\left(\int_{\mathbb{R}^d}K(T-t,z_2)K(T-t,z_1+z_2)dz_2\right)\rho(x-z_1)dz_1,\\
&& I_{\breve\rho}=\mathbb{E}\int_0^t\int_{\mathbb{R}^{2d}}\int_{\mathbb{R}^2}
f_1(x,s,\xi)\bar f_2(y,s,\zeta)(a(\xi)-a(\zeta))\psi(\xi-\zeta)d\xi d\zeta\cdot \nabla_x\breve\rho(x-y)dxdyds,
  \eess
and
   \bess
I_\psi=\frac{1}{2}\int_{\mathbb{R}^{2d}}\breve\rho(x-y)\mathbb{E}\int_0^t\int_{\mathbb{R}^{2}}
\psi(\xi-\zeta)\sum_{k\geq1}|g_k(x,\xi)-g_k(x,\zeta)|^2d\nu_{x,\xi}^1\otimes d\nu_{x,\xi}^2dxdyds.
   \eess
\end{lem}

{\bf Proof.} Since both $\rho$ and $\psi$ have compact support, it is easy to check each
term in (\ref{3.2}) is finite. Set $G_1^2(x,\xi)=\sum_{k=1}^\infty|g_k(x,\xi)|^2$
and $G_2^2(x,\zeta)=\sum_{k=1}^\infty|g_k(x,\zeta)|^2$.


 Let
$\varphi_1(x,\xi)\in C_c^\infty(\mathbb{R}_x^d\times\mathbb{R}_\xi)$ and
$\varphi_2(x,\zeta)\in C_c^\infty(\mathbb{R}_y^d\times\mathbb{R}_\zeta)$.
Recall that
   \bess
&& \tilde\varphi_1(x,\xi)=\int_{\mathbb{R}^d}K(T-x,t-z)\varphi_1(z,\xi)dz,\\
 && \tilde\varphi_2(y,\zeta)=\int_{\mathbb{R}^d}K(T-t,y-z)\varphi_2(z,\zeta)dz.
     \eess
By (\ref{2.10}), we have
    \bess
\langle f_1^+(t),\tilde\varphi_1\rangle=\langle m_1^*,\partial_\xi\tilde\varphi_1\rangle([0,t])+F_1(t),
   \eess
where
    \bess
F_1(t)=\sum_{k\geq1}\int_0^t\int_{\mathbb{R}^{d+1}}g_{k,1}\tilde\varphi_1d\nu_{x,s}^1(\xi)dxd\beta_k(s)
    \eess
and
   \bess
\langle m_1^*,\partial_\xi\tilde\varphi_1\rangle([0,t])&=&\langle f_{1,0},\tilde\varphi_1\rangle
+\int_0^t\langle f_1,a\cdot\nabla_x\tilde\varphi_1\rangle ds\\
&&
+\frac{1}{2}\int_0^t\int_{\mathbb{R}^{d+1}}\partial_\xi\tilde\varphi_1(x,\xi)G_1^2(x,\xi)d\nu^1_{x,s}(\xi)dxds
-\langle m_1,\partial_\xi\tilde\varphi_1\rangle([0,t]).
  \eess

Using Remark \ref{r2.3}, $\langle m_1,\partial_\xi\tilde\varphi_1\rangle(\{0\})=0$
and thus $\langle m_1^*,\partial_\xi\tilde\varphi_1\rangle(\{0\})$ is $\langle f_{1,0},\tilde\varphi_1\rangle$.
Similarly, we have
      \bess
\langle\bar f_2^+(t),\tilde\varphi_2\rangle=\langle\bar m_2^*,\partial_\zeta\tilde\varphi_2\rangle([0,t])+\bar F_2(t),
   \eess
where
    \bess
\bar F_2(t)=-\sum_{k\geq1}\int_0^t\int_{\mathbb{R}^{d+1}}g_{k,2}\tilde\varphi_2d\nu_{y,s}^2(\zeta)dyd\beta_k(s)
    \eess
and
   \bess
\langle\bar m_2^*,\partial_\zeta\tilde\varphi_2\rangle([0,t])&=&\langle\bar f_{2,0},\tilde\varphi_2\rangle
+\int_0^t\langle \bar f_2,a\cdot\nabla_y\tilde\varphi_2\rangle ds\\
&&
-\frac{1}{2}\int_0^t\int_{\mathbb{R}^{d+1}}\partial_\zeta\tilde\varphi_2(y,\zeta)G_2^2(y,\zeta)d\nu^2_{y,s}(\zeta)dyds
+\langle m_2,\partial_\zeta\tilde\varphi_2\rangle([0,t]),
  \eess
where $\langle\bar m_2^*,\partial_\xi\tilde\varphi_2\rangle(\{0\})=\langle\bar f_{2,0},\tilde\varphi_2\rangle$.
Integrating by parts for functions of finite variation
$$\langle m_1^*,\partial_\xi\tilde\varphi_1\rangle([0,t])\langle\bar m_2^*,\partial_\zeta\tilde\varphi_2\rangle([0,t]),$$ we get
   \bess
\langle m_1^*,\partial_\xi\tilde\varphi_1\rangle([0,t])\langle\bar m_2^*,\partial_\zeta\tilde\varphi_2\rangle([0,t])
&=&\langle m_1^*,\partial_\xi\tilde\varphi_1\rangle(\{0\})\langle\bar m_2^*,\partial_\zeta\tilde\varphi_2\rangle(\{0\})\\
&&
+\int_{(0,t]}\langle m_1^*,\partial_\xi\tilde\varphi_1\rangle([0,s))d\langle\bar m_2^*,\partial_\zeta\tilde\varphi_2\rangle(s)\\
&&+\int_{(0,t]}\langle\bar m_2^*,\partial_\zeta\tilde\varphi_2\rangle([0,s))d\langle m_1^*,\partial_\xi\tilde\varphi_1\rangle(s).
   \eess
Since $\bar F_2$ is continuous and $\bar F_2(0)=0$, we have
   \bess
\langle m_1^*,\partial_\xi\tilde\varphi_1\rangle([0,t])\bar F_2(t)=
\int_0^t\langle m_1^*,\partial_\xi\tilde\varphi_1\rangle([0,s])d\bar F_2(s)
+\int_0^t\bar F_2(s)\langle m_1^*,\partial_\xi\tilde\varphi_1\rangle(ds).
   \eess
Denote
$\tilde\sigma=\tilde\varphi_1\tilde\varphi_2 $. Using It\^{o} formula for
$F_1(t)\bar F_2(t)$, we obtain that
   \bess
\langle f_1^+(t),\tilde\varphi_1\rangle\langle\bar f_2^+(t),\tilde\varphi_2\rangle
=\langle\langle f_1^+(t)\bar f_2^+(t),\tilde\sigma \rangle\rangle,
  \eess
where
    \bes
&&\mathbb{E}\langle\langle f_1^+(t)\bar f_2^+(t),\tilde\sigma \rangle\rangle
-\langle\langle f_{1,0}^+\bar f_{2,0}^+,\tilde\sigma_0 \rangle\rangle\nm\\
&=&\mathbb{E}\int_0^t\int_{\mathbb{R}^{2d}}\int_{\mathbb{R}^2}f_1\bar f_2
[a(\xi)\cdot\nabla_x+a(\zeta)\cdot\nabla_y]\tilde\sigma d\xi d\zeta dxdyds\nm\\
&&\mathbb{E}\int_0^t\int_{\mathbb{R}^{2d}}\int_{\mathbb{R}^2}f_1\bar f_2
[(-\Delta_x)^{\frac{\alpha}{2}}+(-\Delta_y)^{\frac{\alpha}{2}}]\tilde\sigma d\xi d\zeta dxdyds\nm\\
&&+\frac{1}{2}\mathbb{E}\int_0^t\int_{\mathbb{R}^{2d}}\int_{\mathbb{R}^2}\partial_\xi\tilde\sigma\bar
f_2(s)G_1^2d\nu_{(x,s)}^1(\xi)d\zeta dxdyds\nm\\
&&-\frac{1}{2}\mathbb{E}\int_0^t\int_{\mathbb{R}^{2d}}\int_{\mathbb{R}^2}\partial_\zeta\tilde\sigma
f_1(s)G_2^2d\nu_{(y,s)}^2(\zeta)d\xi dydxds\nm\\
&&-\mathbb{E}\int_0^t\int_{\mathbb{R}^{2d}}\int_{\mathbb{R}^2}G_{1,2}\tilde\sigma d\nu_{(x,s)}^1(\xi)d\nu_{(y,s)}^2(\zeta)dxdy\nm\\
&&-\mathbb{E}\int_{(0,t]}\int_{\mathbb{R}^{2d}}\int_{\mathbb{R}^2}\bar f^+_2
\partial_\xi\tilde\sigma dm_1(x,s,\xi)d\zeta dy\nm\\
&&+\mathbb{E}\int_{(0,t]}\int_{\mathbb{R}^{2d}}\int_{\mathbb{R}^2} f^-_1\partial_\zeta\tilde\sigma dm_2(y,s,\zeta)d\xi dx,
   \lbl{3.3}
  \ees
where $G_{1,2}(x,y;\xi,\zeta):=\sum_{k\geq 1}g_{k,1}(x,\xi)g_{k,2}(y,\zeta)$ and $\langle\langle\cdot,\cdot\rangle\rangle$
denotes the duality distribution over $\mathbb{R}_x^{d}\times\mathbb{R}_\xi\times\mathbb{R}_y^{d}\times\mathbb{R}$,
and $\tilde\sigma_0 =\tilde\sigma|_{t=0}$.
By a density argument, (\ref{3.3}) remains true for any test function
$\sigma\in C_c^\infty(\mathbb{R}_x^{d}\times\mathbb{R}_\xi\times\mathbb{R}_y^{d}\times\mathbb{R})$.
Let $\sigma=\rho\psi$, where $\rho=\rho(x-y)$ and $\psi=\psi(\xi-\zeta)$.
By the definition of $\tilde \varphi_i$, $i=1,2$, we have
    \bess
\tilde\sigma&=&\int_{\mathbb{R}^d}K(T-t,z_2)\varphi_2(y-z_2,\zeta)dz_2\int_{\mathbb{R}^d}K(T-t,z_1)\varphi_1(x-z_1,\xi)dz_1\\
&=&\int_{\mathbb{R}^d}K(T-t,z_2)\int_{\mathbb{R}^d}K(T-t,z_1)\varphi_1(x-z_1,\xi)\varphi_2(y-z_2,\zeta)dz_1dz_2\\
&=&\psi(\xi-\zeta)\int_{\mathbb{R}^d}K(T-t,z_2)\int_{\mathbb{R}^d}K(T-t,z_1)\rho(x-y-z_1+z_2)dz_1dz_2\\
&=&\psi(\xi-\zeta)\int_{\mathbb{R}^d}\left(\int_{\mathbb{R}^d}K(T-t,z_2)K(T-t,z_1+z_2)dz_2\right)\rho(x-y-z_1)dz_1\\
&=:&\psi(\xi-\zeta)\breve\rho(x-y).
   \eess
Noting that
   \bess
(\nabla_x+\nabla_y)\tilde\sigma=0,\ \ \ \ (\partial_\xi+\partial_\zeta)\tilde\sigma=0,
   \eess
we have
   \bess
&&\mathbb{E}\int_0^t\int_{\mathbb{R}^{2d}}\int_{\mathbb{R}^2}f_1\bar f_2
[a(\xi)\cdot\nabla_x+a(\zeta)\cdot\nabla_y]\tilde\sigma d\xi d\zeta dxdyds\nm\\
&=&\mathbb{E}\int_0^t\int_{\mathbb{R}^{2d}}\int_{\mathbb{R}^2}
f_1(x,s,\xi)\bar f_2(y,s,\zeta)(a(\xi)-a(\zeta))\psi(\xi-\zeta)d\xi d\zeta\cdot \nabla_x\breve\rho(x-y)dxdyds.
   \eess
The last term in (\ref{3.3}) is
   \bess
&&\mathbb{E}\int_{(0,t]}\int_{\mathbb{R}^{2d}}\int_{\mathbb{R}^2} f^-_1\partial_\zeta\tilde\sigma dm_2(y,s,\zeta)d\xi dx
\\&=&-\mathbb{E}\int_{(0,t]}\int_{\mathbb{R}^{2d}}\int_{\mathbb{R}^2} f^-_1\partial_\xi\tilde\sigma dm_2(y,s,\zeta)d\xi dx\\
&=&-\mathbb{E}\int_{(0,t]}\int_{\mathbb{R}^{2d}}\int_{\mathbb{R}^2} \tilde\sigma d\nu_{(x,s)}^{1,-} dm_2(y,s,\zeta)d\xi dx\\&\leq&0
   \eess
since $\alpha\geq0$ and $m_2,\,n_2$ are non-negative measure. Similarly, we have
   \bess
-\mathbb{E}\int_{(0,t]}\int_{\mathbb{R}^{2d}}\int_{\mathbb{R}^2}\bar f^+_2\partial_\xi\tilde\sigma dm_1(x,s,\xi)d\zeta dy=-\mathbb{E}\int_{(0,t]}\int_{\mathbb{R}^{2d}}\int_{\mathbb{R}^2}\tilde\sigma d\nu_{(y,s)}^{2,+} dm_1(x,s,\xi)d\zeta dy\leq0.
  \eess
Integrating by part, we get
   \bess
&&\frac{1}{2}\mathbb{E}\int_0^t\int_{\mathbb{R}^{2d}}\int_{\mathbb{R}^2}\partial_\xi\tilde\sigma\bar
f_2(s)G_1^2d\nu_{(x,s)}^1(\xi)d\zeta dxdyds\nm\\
&&-\frac{1}{2}\mathbb{E}\int_0^t\int_{\mathbb{R}^{2d}}\int_{\mathbb{R}^2}\partial_\zeta\tilde\sigma
f_1(s)G_2^2d\nu_{(y,s)}^2(\zeta)d\xi dydxds\nm\\
&&-\mathbb{E}\int_0^t\int_{\mathbb{R}^{2d}}\int_{\mathbb{R}^2}G_{1,2}\tilde\sigma d\nu_{(x,s)}^1(\xi)d\nu_{(y,s)}^2(\zeta)dxdy\nm\\
&=&\frac{1}{2}\mathbb{E}\int_0^t\int_{\mathbb{R}^{2d}}\int_{\mathbb{R}^2}\tilde\sigma(G_1^2-2G_{1,2}+G_2^2)
d\nu_{(x,s)}^1\otimes d\nu_{(y,s)}^2(\xi,\zeta)dxdyds\\
&=&\frac{1}{2}\int_{\mathbb{R}^{2d}}\breve\rho(x-y)\mathbb{E}\int_0^t\int_{\mathbb{R}^{2}}
\psi(\xi-\zeta)\\
&&\qquad\qquad\qquad\times\sum_{k\geq1}|g_k(x,\xi)-g_k(x,\zeta)|^2d\nu_{x,\xi}^1\otimes d\nu_{x,\xi}^2dxdyds\\
&=:&I_\psi.
   \eess
Combining the above discussion, we obtain the desired results. $\Box$

{\bf Proof of Theorem \ref{t3.1}} We first use the Lemma \ref{l3.1} to prove the uniqueness.
The additive case: $\Phi(u)$ independent on $u$. Let $f_i,\,i=1,2$ be two generalized solution
to (\ref{1.1}). Then, we use (\ref{3.2}) with $g_k$ independent of $\xi$ and $\zeta$. By (\ref{1.3})
the last term $I_\psi$ is bounded by
    \bess
\frac{tD_1}{2}\|\psi\|_{L^\infty}\int_{\mathbb{R}^{2d}}|x-y|^2\breve\rho(x-y)dxdy.
   \eess
Note that if $\rho(x)\equiv C$, by using the properties of the heat kernel $K$, we have
   \bess
\breve\rho(x)&=&\int_{\mathbb{R}^d}\left(\int_{\mathbb{R}^d}K(T-t,z_2)K(T-t,z_1+z_2)dz_2\right)\rho(x-z_1)dz_1\\
&=&C\int_{\mathbb{R}^d}K(T-t,z_2)\left(\int_{\mathbb{R}^d}K(T-t,z_1+z_2)dz_1\right)dz_2\\
&=&C.
   \eess
Taking $\psi:=\psi_\delta$ and $\rho=\rho_\epsilon$, where $(\psi_\delta)$ and
$(\rho_\epsilon)$ are approximations to the identity on $\mathbb{R}$ and $\mathbb{R}^d$ respectively,
we obtain
   \bes
I_\psi\leq\frac{t\tilde D_1}{2}\epsilon^2\delta^{-1},
   \lbl{3.4}\ees
where $\tilde D_1=D_1\int_{\mathbb{R}^d}z^2\rho(z)dz<\infty$ because of the compact support of $\rho$.
Let $t\in[0,T]$, $t_n\downarrow t$ and $\nu_{x,t}^{i,+}$ be a weak-limit of $\nu_{x,t_n}^{i,+}$ in
sense of (\ref{2.6}). Then
$\nu_{x,t}^{i,+}$ satisfies
   \bess
\mathbb{E}\int_{\mathbb{R}^d}\int_{\mathbb{R}}|\xi|^pd\nu_{x,t}^{i,+}dx\leq C_p,
   \eess
and we have a similar bound for $\nu^{i,-}$. Denote
    \bess
\tilde K(T-x,t)=\int_{\mathbb{R}^d}K(T-t,z)K(T-x,t+z)dz.
   \eess
We can rewrite the integration as
   \bess
&&\mathbb{E}\int_{\mathbb{R}^{d}}\int_{\mathbb{R}^{d}}\int_{\mathbb{R}}\tilde K(T-x,t-y)f_1^{\pm}(x,t,\xi)
\bar f_2^\pm(x,t,\xi)d\xi  dydx \\
&=&\mathbb{E}\int_{\mathbb{R}^{2d}}\int_{\mathbb{R}^2}\breve\rho_\epsilon(x-y)\psi_\delta(\xi-\zeta)f_1^{\pm}(x,t,\xi)
\bar f_2^\pm(y,t,\zeta)d\xi d\zeta dxdy+\eta_t(\epsilon,\delta),
   \eess
where $\lim\limits_{\epsilon,\delta\to0}\eta_t(\epsilon,\delta)=0$. Now, we need a bound on the
term $I_\rho$. Since $a$ has at most polynomial growth, similar to the proof of \cite[Theorem 11, pp 1029]{DV},
there exists a positive constant $C_p$ such that
    \bes
&&\Big|\mathbb{E}\int_0^t\int_{\mathbb{R}^{2d}}\int_{\mathbb{R}^2}
f_1(x,s,\xi)\bar f_2(y,s,\zeta)(a(\xi)-a(\zeta))\psi_\delta(\xi-\zeta)d\xi d\zeta\nm\\[1.5mm]
&&\qquad\qquad\qquad\cdot \nabla_x\breve\rho_\epsilon(x-y)dxdyds\Big|
\leq tC_p\delta\epsilon^{-1}.
   \lbl{3.5}\ees
We then gather (\ref{3.4}), (\ref{3.5}) and (\ref{3.2}) to deduce for $t\in[0,T]$
   \bess
&&\mathbb{E}\int_{\mathbb{R}^{d}}\int_{\mathbb{R}^{d}}\int_{\mathbb{R}}\tilde K(T-x,t-y)f_1^{\pm}(x,t,\xi)
\bar f_2^\pm(x,t,\xi)d\xi  dydx\\
&\leq&\mathbb{E}\int_{\mathbb{R}^{d}}\int_{\mathbb{R}^{d}}\int_{\mathbb{R}}\tilde K(T,x-y)f_{1,0}
\bar f_{2,0}d\xi  dydx+r(\epsilon,\delta),
   \eess
where the remainder $r(\epsilon,\delta)$ is
    \bess
 r(\epsilon,\delta)=TC_p\delta\epsilon^{-1}+\frac{T\tilde D_1}{2}\epsilon^2\delta^{-1}+
 \eta_t(\epsilon,\delta)+\eta_0(\epsilon,\delta).
    \eess
Taking $\delta=\epsilon^{\frac{4}{3}}$ and letting $\epsilon\to0$, we have
$\lim\limits_{\epsilon\to0}r(\epsilon,\delta)=0$ and
   \bes
&&\mathbb{E}\int_{\mathbb{R}^{d}}\int_{\mathbb{R}^{d}}\int_{\mathbb{R}}\tilde K(T-x,t-y)f_1^{\pm}(x,t,\xi)
\bar f_2^\pm(x,t,\xi)d\xi  dydx\nm\\
&\leq&\mathbb{E}\int_{\mathbb{R}^{d}}\int_{\mathbb{R}^{d}}\int_{\mathbb{R}}\tilde K(T,x-y)f_{1,0}
\bar f_{2,0}d\xi  dydx.
   \lbl{3.6}\ees
Assume that $f$ is a generalized solution to (\ref{1.1}) with initial datum $\textbf{1}_{u_0>\xi}$.
Since $f_0$ is the Heaviside function, we get the identity $f_0\bar f_0=0$. Taking
$f_1=f_2=f$, by the positive property of $K$ ($K(x,t)>0$ for any $t>0$ and $x\in\mathbb{R}^d$),
we deduce that $f^+(1-f^+)=0$ a.e., i.e. $f^+\in\{0,1\}$ a.e.. The fact $-\partial_\xi f^+$
is a Young measure gives the conclusion: indeed, by Fubini Theorem, for any $t\in[0,T]$, there
is a set $E_t$ of full measure in $\mathbb{R}^d\times\Omega$ such that, for $(x,\omega)\in E_t$,
$f^+(x,t,\xi,\omega)\in\{0,1\}$ for a.e. $\xi\in\mathbb{R}$. Recall that $-\partial_\xi f^+(x,t,\cdot,\omega)$
is a probability measure on $\mathbb{R}$ so that there exists $u^+(x,t,\omega)\in\mathbb{R}$ such that
$f^+(x,t,\xi,\omega)=\textbf{1}_{u^+(x,t,\omega)>\xi}$ for almost every $(x,\xi,\omega)$. In particular,
$u^+=\int_\mathbb{R}(f^+-\textbf{1}_{\xi>0})d\xi$ for almost every $(x,\omega)$. A similar result also
holds for $f^-$.

It follows from the discussion after Definition \ref{d2.6} that
$f^+$ being solution in the sense of Definition \ref{d2.6} implies that
$u^+$ is a solution in the sense of Definition \ref{d2.3}. Since $f=f^+$ a.e.,
this shows the reduction of reduction of generalized solutions to solutions. If now
$u_1$ and $u_2$ are two solutions to (\ref{1.1}), we deduce from (\ref{3.6}) with
$f_i=\textbf{1}_{u_i>\xi}$ and from the identity
   \bess
\int_\mathbb{R}\textbf{1}_{u_1>\xi}\overline{\textbf{1}_{u_2>\xi}}d\xi=(u_1-u_2)^+
   \eess
the contraction property
   \bess
&&\mathbb{E}\int_{\mathbb{R}^d}\left[\tilde K(T-t,\cdot)\ast(u_1(t,\cdot)-u_2(t,\cdot))^+\right](x)dx\\
&\leq& \mathbb{E}\int_{\mathbb{R}^d}\left[\tilde K(T,\cdot)\ast(u_1(0,\cdot)-u_2(0,\cdot))^+\right](x)dx,
    \eess
which implies the uniqueness of solutions. Actually, due to $K(T-x,t)>0$ for any $t\in[0,T]$ and every
$x\in\mathbb{R}^d$, we have if $u_{1,0}=u_{2,0}$, then $u_1=u_2$.

In the multiplicative case ($\Phi$ depending on $u$), the reasoning is similar, except
that there is an additional term in the bound on $I_\psi$. More precisely, by Hypothesis
(\ref{1.4}) we obtain in place of (\ref{3.4}) the estimate
   \bess
I_\psi\leq\frac{T\tilde D_1}{2}\epsilon^2\delta^{-1}+\frac{D_1}{2}I_\psi^h,
   \eess
where
   \bess
I_\psi^h=\mathbb{E}\int_0^t\int_{\mathbb{R}^{2d}}\breve\rho_\epsilon\int_{\mathbb{R}^2}
\psi_\delta(\xi-\zeta)|\xi-\zeta|h(|\xi-\zeta|)d\nu_{x,\sigma}^1\otimes \nu_{y,\sigma}^2(\xi,\zeta)dxdyd\sigma.
   \eess
Choosing $\psi_\delta(\xi)=\delta^{-1}\psi_1(\delta^{-1}\xi)$ with $\psi_1$ compactly supported gives
   \bess
I_\psi\leq\frac{T\tilde D_1}{2}\epsilon^2\delta^{-1}+\frac{TD_1C_\psi h(\delta)}{2},\ \
C_\psi:=\sup_{\xi\in\mathbb{R}}\|\xi\psi_1(\xi)\|,
   \eess
which implies that  $\lim\limits_{\epsilon\to0,\delta=\epsilon^{4/3}}I_\psi=0$
Similar to the additional case and the proof of Theorem 11 in \cite{DV}, one can finish the proof of
uniqueness of solution, which is the part
of Theorem \ref{t3.1}.

\textbf{(Existence)} We prove the existence by a vanishing viscosity method.  Assume that  $u_0\in L^\infty\cap L^1\cap BV(\mR^d)$.
 \vskip1mm\par
Consider the Cauchy problem:
\begin{eqnarray}\lbl{3.7}
\left\{ \begin{array}{lll}
du^\varepsilon(x,t)+[\div_xA(u^\varepsilon)
+(-\Delta_x)^{\frac{\alpha}{2}}u^\varepsilon]dt\\[1.5mm]
\qquad\qquad\qquad\qquad\qquad\qquad-
\varepsilon\Delta u^\varepsilon dt=\Phi^\varepsilon(u^\varepsilon)dW(t), \ \ (x,t)\in (0,T) \times \mR^d,
\\[1.5mm] u^\varepsilon(t=0)=u_0, \ \ x\in\mR^d,
\end{array}\right.
\end{eqnarray}
where $\Phi^\varepsilon$ is a suitable Lipschitz approximation of $\Phi$ satisfying
(\ref{1.3}) and (\ref{1.4}) uniformly. We define $g_k^\varepsilon$ and $G_\varepsilon$
as in the case $\varepsilon=0$.

Similar to the proof of Lemma 4.9 in \cite{FN}, we can prove equation
(\ref{3.7}) has a unique solution
$u^\varepsilon\in L^\infty([0,T],L^p(\mathbb{R}^d)\cap L^2([0,T],H^{\frac{\alpha}{2}}(\mathbb{R}^d))$ provided
that $u_0\in L^p(\mathbb{R}^d)$, $p\geq2$. Moreover, by using It\^{o} formula, one can prove
that $u^\varepsilon$ satisfies the energy inequality
    \bess
&&\mathbb{E}\left[\|u^\varepsilon(t)\|_2^2\right]+2\varepsilon\mathbb{E}\int_0^t\int_{\mathbb{R}^d}|\nabla u^\varepsilon(s,x)|^2dxds\\
&&+2\mathbb{E}\int_0^t\int_{\mathbb{R}^d}u^\varepsilon(s,x)
(-\Delta)^{\frac{\alpha}{2}}u^\varepsilon(s,x)dxds\\
&=&\mathbb{E}\left[\|u(0)\|_2^2\right]-
2\mathbb{E}\int_0^t\int_{\mathbb{R}^d}u^\varepsilon(s,x)\div A(u^\varepsilon(s,x))dxds\\
&&
+\mathbb{E}\int_0^t\int_{\mathbb{R}^d}G_\varepsilon^2(u^\varepsilon(s,x))dxds\\
&\leq&\mathbb{E}\left[\|u(0)\|_2^2+\|\hat g\|_{L^1}\right]
+D_0\int_0^t\mathbb{E}\left[\|u^\varepsilon(s)\|_2^2\right]ds,
   \eess
which implies that by using Gronwall's Lemma
     \bes
&&\mathbb{E}\left[\|u^\varepsilon(t)\|_2^2\right]+
2\varepsilon\mathbb{E}\int_0^t\|\nabla u^\varepsilon(s)\|_{L^2(\mathbb{R}^d)}^2ds\nm\\
&&
+2\mathbb{E}\int_0^t\|u^\varepsilon(s)\|_{H^{\frac{\alpha}{2}}(\mathbb{R}^d)}^2ds
\leq C_T\mathbb{E}\left[\|u(0)\|_2^2+\|\hat g\|_{L^1}\right].
   \lbl{j.1}\ees
Also, for $p\geq2$, by It\^{o} formula applied to $|u^\varepsilon|^p$ and a martingale inequality
    \bes
\mathbb{E}\left(\sup_{t\in[0,T]}\|u^\varepsilon(t)\|_{L^p(\mathbb{R}^d)}^p\right)+
\varepsilon\mathbb{E}\int_0^T\int_{\mathbb{R}^d}|u^\varepsilon(x,t)|^{p-2}|\nabla u^\varepsilon(x,t)|^2dxdt\leq
C(p,u_0,T).
    \lbl{j.2}\ees
Similar to the discussion in Section 2 and the proof of Proposition 18 in \cite{DV}, we can obtain
the following result.
  \begin{prop}\lbl{p3.1} (Kinetic formulation) Let $u_0\in L^\infty\cap L^1\cap BV(\mR^d)$ and
let $u^\varepsilon$ be the solution to (\ref{3.7}). Then $f^\varepsilon:=\textbf{1}_{u^\varepsilon>\xi}$
satisfies: for all $\varphi\in C_c^2(\mathbb{R}^d\times[0,T)\times\mathbb{R})$,
    \bes
&&\int_0^T\langle f^\varepsilon(t),\partial_t\varphi(t)\rangle dt+\langle f_0,\varphi(0)\rangle
+\int_0^T\langle f^\varepsilon(t),a(\xi)\cdot\nabla\varphi(t)+(-\Delta)^{\frac{\alpha}{2}}\varphi(t)-\varepsilon\Delta\varphi(t)\rangle dt\nm\\
&=&-\sum_{k\geq1}\int_0^T\int_{\mathbb{R}^{d+1}}g_k^\varepsilon(x,\xi)\varphi(x,t,\xi)d\nu_{x,t}^\varepsilon(\xi)dxd\beta_k(t)\nm\\
&&-\frac{1}{2}\int_0^T\int_{\mathbb{R}^{d+1}}\partial_\xi\varphi(x,t,\xi)G^2_\varepsilon(x,\xi)d\nu_{x,t}^\varepsilon(\xi)dxdt+
m^\varepsilon(\partial_\xi\varphi),
    \lbl{3.8}\ees
a.s., where $f_0(\xi)=\textbf{1}_{u_0>\xi}$, $m^\varepsilon=m^\varepsilon_1+m^\varepsilon_2$, $m^\varepsilon_1$ is defined as (\ref{2.4}) and
   \bess
\nu_{x,t}^\varepsilon=\delta_{u^\varepsilon(x,t)},\ \ \
m^\varepsilon_2=\varepsilon
|\nabla u^\varepsilon(x,t)|^2\delta_{u^\varepsilon(x,t)=\xi}.
       \eess
\end{prop}

Let $\eta_\epsilon$ satisfy the assumption in Definition \ref{d2.1} and $\eta_\epsilon(r)\to|r|$
as $\epsilon\to0$. It\^{o} formula gives
   \bes
d\eta_\epsilon(u^\varepsilon)&=&-\eta'_\epsilon(u^\varepsilon)[\div_xA(u^\varepsilon)
-(-\Delta_x)^{\frac{\alpha}{2}}u^\varepsilon+
\varepsilon\Delta u^\varepsilon] dt\nm\\
&&+\eta'_\epsilon(u^\varepsilon)\Phi^\varepsilon(u^\varepsilon)dW(t)
+\frac{1}{2}\eta''_\epsilon(u^\varepsilon)G^2_\varepsilon dt
   \lbl{3.9}\ees
The convex of $\eta$ implies that
   \bess
\varepsilon\eta^\prime(u^\varepsilon)\Delta u^\varepsilon&=&\varepsilon
\Delta\eta(u^\varepsilon)-\varepsilon\eta^{\prime\prime}(u^\varepsilon)
|\nabla u^\varepsilon|^2
\leq \Delta\eta(u^\varepsilon)\\[1.5mm]
\eta^\prime(u^\varepsilon)(-\Delta_x)^{\frac{\alpha}{2}}u^\varepsilon(x,t)
&\geq&
c_0\int_{\mR^d}\frac{\eta(u^\varepsilon(x,t))-
\eta(u^\varepsilon(t,z+x))}
{|z|^{d+\alpha}}dz
=(-\Delta_x)^{\frac{\alpha}{2}}\eta(u^\varepsilon(x,t)).
  \eess
Integrating (\ref{3.9}) over $\mathbb{R}^d$, using the above two
inequalities, taking expectation and letting $\epsilon\to0$, we get
    \bess
\mathbb{E}\int_{\mathbb{R}^d}|u^\varepsilon(x,t)|dx&\leq&\mathbb{E}
\int_{\mathbb{R}^d}|u_0(x)|dx
-\mathbb{E}\int_{\mathbb{R}^d}\mbox{sgn}(u^\varepsilon)[\div_xA(u^\varepsilon)] dt\\
&&
+\frac{1}{2}\lim\limits_{\epsilon\to0}\mathbb{E}\int_{\mathbb{R}^d}\eta''_\epsilon(u^\varepsilon)G^2_\varepsilon dt\\
&\leq&\|u_0\|_{L^1(\mathbb{R}^d)}+\frac{D_0}{2}\|\hat g\|_{L^1(\mathbb{R}^d)}.
   \eess
which implies that $u^\varepsilon\in L^1(\mathbb{R}^d)$.

It follows from (\ref{j.1}) that $u^\varepsilon$ weakly converges in $H^{\frac{\alpha}{2}}(\mathbb{R}^d)$.

Equation (\ref{3.9}) is close to the kinetic equation (\ref{2.8}) satisfied by the solution
to (\ref{1.1}). For $\varepsilon\to0$, we lose the precise structures of
$m^\varepsilon=\varepsilon|\nabla u^\varepsilon|^2\delta_{u^\varepsilon=\xi}$ and $n^\varepsilon$, and
obtain a solution $u$ to (\ref{1.1}). More precisely, we will prove the following result.
\begin{theo}\lbl{t3.2} (Convergence of the parabolic approximation) Let $u_0\in L^\infty\cap L^1(\mathbb{R}^d)$.
There exists a unique solution $u$ to (\ref{1.1}) with initial datum $u_0$ which
is the strong limit of $(u^\varepsilon)$ as $\varepsilon\to0$: for every $T>0$, for every
$1\leq p<\infty$,
   \bess
\lim\limits_{\varepsilon\to0}\mathbb{E}\|u^\varepsilon-u\|_{L^p(\mathbb{R}^d\times(0,T))}=0.
   \eess
 Moreover, $(u^\varepsilon)$ converges weakly to $u$ in $H^{\frac{\alpha}{2}}(\mathbb{R}^d)$.
  \end{theo}
The proof of Theorem \ref{t3.2} is a straightforward consequence of both the result of
reduction of generalized solution to solution (uniqueness of Theorem \ref{t3.1}) and a priori
estimates derived in the following.

{\bf Estimates of $m^\varepsilon_1$ and $m^\varepsilon_2$:} similar to that in \cite{DV},
we analyze $m^\varepsilon_1$ and $m^\varepsilon_2$.
By (\ref{j.1}), we obtain a uniform bound $\mathbb{E}m^\varepsilon_2(\mathbb{R}^d\times[0,T]\times\mathbb{R})\leq C$.
Furthermore, the second term in the left hand-side of (\ref{j.2}) is
$\mathbb{E}\int_{\mathbb{R}^d\times[0,T]\times\mathbb{R}}|\xi|^{p-2}dm_2^\varepsilon(x,t,\xi)$,
so we have
   \bess
\mathbb{E}\int_{\mathbb{R}^d\times[0,T]\times\mathbb{R}}|\xi|^{p}dm_2^\varepsilon(x,t,\xi)\leq C_p.
   \eess
We also have the the improved estimate, for $p\geq0$
   \bes
\mathbb{E}\Big|\int_{\mathbb{R}^d\times[0,T]\times\mathbb{R}}|\xi|^{2p}d
m^\varepsilon(x,t,\xi)\Big|^2
\leq C_p,
   \lbl{3.12}\ees
where for $\psi\in C^2$ and $\psi''\geq0$,
   \bess
\int_{\mathbb{R}^d\times[0,T]\times\mathbb{R}}\psi(\xi)d
m_1^\varepsilon(x,t,\xi)=\int_{\mathbb{R}^d\times[0,T]}\left(\psi'(u^\varepsilon)(-\Delta)^{\frac{\alpha}{2}}u^\varepsilon
-(-\Delta)^{\frac{\alpha}{2}}\psi(u^\varepsilon)\right)dxdt.
   \eess
To prove (\ref{3.12}), we apply It\^{o} formula to $\psi(u^\varepsilon)$, $\psi(\xi):=|\xi|^{2p+2}$
   \bess
&&d\psi(u^\varepsilon)+\div (U)dt+\varepsilon\psi''(u^\varepsilon)|\nabla u^\varepsilon|^2dt\\
&&+[\psi'(u^\varepsilon)(-\Delta)^{\frac{\alpha}{2}}u^\varepsilon
-(-\Delta)^{\frac{\alpha}{2}}\psi(u^\varepsilon)]dt\\
&=&\psi'(u^\varepsilon)\Phi_\varepsilon(u^\varepsilon)dW+\frac{1}{2}\psi''(u^\varepsilon)G_\varepsilon^2dt
-(-\Delta)^{\frac{\alpha}{2}}\psi(u^\varepsilon)dt,
   \eess
where
$U=:\int_0^{u^\varepsilon}a^\varepsilon(\xi)\psi'(\xi)d\xi-\varepsilon\nabla\psi(u^\varepsilon)$.
Integrating over $\mathbb{R}^d\times[0,T]$ yields that
   \bess
&&\int_0^T\int_{\mathbb{R}^d} [\varepsilon\psi''(u^\varepsilon)|\nabla u^\varepsilon|^2
+\psi'(u^\varepsilon)(-\Delta)^{\frac{\alpha}{2}}u^\varepsilon
-(-\Delta)^{\frac{\alpha}{2}}\psi(u^\varepsilon)]dxdt\\
&\leq&\int_{\mathbb{R}^d}\psi(u_0)dx
+\sum_{k\geq1}\int_0^T\int_{\mathbb{R}^d}\psi'(u^\varepsilon)g_{k,\varepsilon}(x,u^\varepsilon)dxd\beta_x(t)\\
&&+\frac{1}{2}\int_0^T\int_{\mathbb{R}^d}\psi''(u^\varepsilon)G^2_{\varepsilon}(x,u^\varepsilon)dxdt.
   \eess
Taking the square, then expectation, we deduce by It\^{o} isometry
   \bess
&&\mathbb{E}\Big|\int_0^T\int_{\mathbb{R}^d} [\varepsilon\psi''(u^\varepsilon)|\nabla u^\varepsilon|^2
+\psi'(u^\varepsilon)(-\Delta)^{\frac{\alpha}{2}}u^\varepsilon
-(-\Delta)^{\frac{\alpha}{2}}\psi(u^\varepsilon)]dxdt\Big|^2\\
&\leq&3\mathbb{E}\Big|\int_{\mathbb{R}^d}\psi(u_0)dx\Big|^2
+3\mathbb{E}\int_0^T\sum_{k\geq1}\Big|\int_{\mathbb{R}^d}\psi'(u^\varepsilon)g_{k,\varepsilon}(x,u^\varepsilon)dx\Big|^2dt\\
&&+\frac{3}{2}\mathbb{E}\Big|\int_0^T\int_{\mathbb{R}^d}\psi''(u^\varepsilon)G^2_{\varepsilon}(x,u^\varepsilon)dxdt\Big|^2.
   \eess
By using Cauchy-Schwarz inequality and the condition (\ref{1.3}), we obtain (\ref{3.12}).

{\bf Estimate on $\nu^\varepsilon$:} This part is similar to that in \cite{DV}. Using the bound (\ref{j.2}),
we get
   \bes
\mathbb{E}\sup_{t\in[0,T]}\int_{\mathbb{R}^d}\int_{\mathbb{R}}|\xi|^pd\nu^\varepsilon_{x,t}(\xi)dx\leq C_p
   \lbl{3.13}\ees
and, in particular,
   \bes
\mathbb{E}\int_0^T\int_{\mathbb{R}^d}\int_{\mathbb{R}}|\xi|^pd\nu^\varepsilon_{x,t}(\xi)dx\leq C_p.
   \lbl{3.14}\ees

Consider a sequence $(\varepsilon_n)\downarrow0$. First, by (\ref{3.14}) and the proposition
\ref{p2.1} (see \cite[Theorem 5 and Corollary 6]{DV}), the convergence $\nu^{\varepsilon_n}\to\nu$
and $f^{\varepsilon_n}\rightharpoonup f$ in $L^\infty(\Omega\times\mathbb{R}^d\times(0,T)\times\mathbb{R})$-weak$-^*$.
Besides, the bound (\ref{3.13}) is stable: $\nu$ satisfies (\ref{2.7}).

For $r\in\mathbb{N}^*$, let $K_r=\mathbb{R}^d\times[0,T]\times[-r,r]$ and let $\mathcal {M}_r$ denote
the space of bounded Borel measures over $K_r$ (with norm given by the total variation of measures).
It is the topological dual of $C(K_r)$, the set of continuous functions on $K_p$. Since
$\mathcal {M}_r$ is separable, the space $L^2(\Omega;\mathcal {M}_r)$ is the topological dual
space of $L^2(\Omega,C(K_r))$, see Th$\acute{e}$or$\grave{e}$m 1.4.1 in \cite{Dro01}.
The estimate (\ref{3.12}) with $p=0$ gives a uniform bound on $m^{\varepsilon_n}$ and $n^{\varepsilon_n}$
in $L^2(\Omega;\mathcal {M}_r)$: there exists $m_r\in L^2(\Omega;\mathcal {M}_r)$ such
that up to subsequence, $m^{\varepsilon_n}\rightharpoonup m$
in $L^2(\Omega;\mathcal {M}_r)$-weak star. By a diagonal process, we obtain
$m_r=m_{r+1}$ in $L^2(\Omega;\mathcal {M}_r)$ and the convergence in
all the spaces $L^2(\Omega;\mathcal {M}_r)$-weak star of a single subsequence still
denoted $(m^{\varepsilon_n})$. The condition at infinity
(\ref{3.12}) shows that $m$ defines two elements of $L^2(\Omega;\mathcal {M})$,
where $\mathcal {M}$ denotes the space of bounded Borel measures over $\mathbb{R}^d\times(0,T)\times\mathbb{R}$.
It follows that
   \bes
\mathbb{E}\Big|\int_{\mathbb{R}^d\times[0,T]\times\mathbb{R}}|\xi|^{2p}dm^\varepsilon(x,t,\xi)
\Big|^2
\leq C_p,
   \lbl{3.15}\ees
which is exactly as (45) in \cite{DV}. So following the idea of
\cite{DV}, we can prove the measure $m$ satisfies $1,\,2,\,$ and $3$ in Definition \ref{d2.2},
that is, $m$ is a kinetic measure.

{\bf Proof of Theorem \ref{t3.2}} By the proof of uniqueness, there corresponds a solution $u$
to this $f:f=\textbf{1}_{u>\xi}$. This proves the existence of a solution $u$ to (\ref{1.1}). Besides,
owing to the particular structure of $f^\varepsilon$ and $f$, we have
   \bess
\|u^{\varepsilon_n}\|_{L^2(\mathbb{R}^d\times(0,T))}^2- \|u\|_{L^2(\mathbb{R}^d\times(0,T))}^2
=\int_0^T\int_{\mathbb{R}^{d+1}}2\xi(f^{\varepsilon_n}-f)d\xi dxdt
   \eess
and (using the bound on $u^\varepsilon$ in $L^3(\mathbb{R}^d)$)
   \bess
\mathbb{E}\int_0^T\int_{\mathbb{R}^{d}}\int_{|\xi|>R}|2\xi(f^{\varepsilon_n}-f)|d\xi dxdt\leq
\frac{C}{1+R}.
   \eess
It follows that $u^{\varepsilon_n}$ converges in norm to $u$ in the Hilbert space
$L^2(\Omega\times\mathbb{R}^d\times(0,T))$. Using the weak convergence, we deduce the
strong convergence. Since $u$ is unique, the whole sequence actually converges. This
gives the result of theorem for $p=2$. The case of general $p$ follows from the bound
on $u^\varepsilon$ in $L^q$ for arbitrary $q$  and H\"{o}lder inequality.
Moreover, it follows from the uniform bound (\ref{j.1}) that
$(u^\varepsilon)$ converges weakly to $u$ in $H^{\frac{\alpha}{2}}(\mathbb{R}^d)$. This completes the
proof of Theorem \ref{t3.2}. $\Box$

The existence of solution in sense of Definition \ref{d2.2} are obtained by Theorem \ref{t3.2}.
This completes the proof of Theorem \ref{t3.1}.  $\Box$

\vskip4mm\noindent
\textbf{\large{Acknowledgements}}
 \vskip3mm\par
This research was partly supported by the NSF of China grants 11501577, 11301146 and 11531006.

\end{document}